\magnification=1200
\font\fp=cmr9
\parskip=10pt plus 1 pt
\overfullrule=0pt
\input amssym

%%%%%%% EXPRESIONES MATEMATICAS %%%%%%%

%%% funciones

%%% sucesiones

\def\sn#1#2{( #1_n ) {\atop {n \geq #2}}}

%%% productos

\def\lp{\big< \kern -.2 cm \big<}
\def\rp{\big> \kern -.2 cm \big>}

%%% demostraciones

\def\qed{\ifhmode\unskip\nobreak\fi\ifmmode\ifinner\else\hskip5 pt \fi\fi
\hbox{\hskip5 pt \vrule width4 pt  height6 pt  depth1.5 pt \hskip 1 pt }}

%%%%%%% ESTILOS DE LETRA %%%%%%%

%%% negrita

\def\b#1{{\bf #1}}

%%% funciones

%\font\Func=Chancery at 12 bp
%\def\F#1{\hbox{\Func #1}}
%\font\func=Chancery at 6 bp
%\def\f#1{\hbox{\func #1}}

%%% conjuntos

\font\msbmnormal=msbm10
\font\msbmpeq=msbm7
\font\msbmmuypeq=msbm5
\newfam\numeros
\textfont\numeros=\msbmnormal
\scriptfont\numeros=\msbmpeq
\scriptscriptfont\numeros=\msbmmuypeq
\def\num{\fam\numeros\msbmnormal}

\def\R{{\num R}}
\def\C{{\num C}}

\def\D{{\num D}}
\def\T{{\num T}}

%%%%%%%% O R L A S %%%%%%% EJEMPLO: $\orlar{4563 \int Ml}$

\def\boxit#1#2{\setbox1=\hbox{\kern#1{#2}\kern#1}
\dimen1=\ht1 \advance\dimen1 by #1 \dimen2=\dp1 \advance\dimen2 by #1
\setbox1=\hbox{\vrule height\dimen1 depth\dimen2\box1\vrule}%
\setbox1=\vbox{\hrule\box1\hrule}%
\advance\dimen1 by .4pt \ht1=\dimen1
\advance\dimen2 by .4pt \dp1=\dimen2 \box1\relax}
\def\orlar#1{\ifmmode\boxit{3pt}{$#1$}
\else
\boxit{3pt}{#1}\fi}

%%%%%%%%%%%%%%%%%%%%%%%%%%%%%%%%%%%%%%%%%%%%%%%%%%%%%%%%%%%%%%%%%%%%%%%%%%%%%%%%%%
%%%%%%%%%%%%%%%%%%%%%%%%%%%%%%   AQUI EMPIEZA EL TEXTO   %%%%%%%%%%%%%%%%%%%%%%%%%
%%%%%%%%%%%%%%%%%%%%%%%%%%%%%%%%%%%%%%%%%%%%%%%%%%%%%%%%%%%%%%%%%%%%%%%%%%%%%%%%%%

\hsize 15true cm
\centerline{\bf A connection between orthogonal polynomials on the unit circle}
\vskip 0.2cm
\centerline{\bf and matrix orthogonal polynomials on the real line}

$$
\vbox{\halign {&\quad \hfil # \hfil \quad \cr
{\it M.J. Cantero},$\;\;${\it M.P. Ferrer},$\;\;$
{\it L. Moral},$\;\;${\it L. Vel\'azquez}
\cr \cr
{\fp Departamento de Matem\'atica Aplicada. Universidad de Zaragoza. Spain.}
\cr}}
$$

\vskip 0.25cm

\noindent {\bf Abstract}

{\fp Szeg\H o's procedure to connect orthogonal polynomials on the unit circle
and or\-tho\-go\-nal polynomials on $[-1,1]$ is generalized to nonsymmetric measures.
It generates the so-called semi-orthogonal functions on the linear space of Laurent
polynomials $\Lambda$, and leads to a new orthogonality structure in the module
$\Lambda \times \Lambda$. This structure can be interpreted in terms of a
$2 \times 2$ matrix measure on $[-1,1]$, and semi-orthogonal functions provide the
corresponding sequence of orthogonal matrix polynomials. This gives a connection
between orthogonal polynomials on the unit circle and certain classes of matrix
orthogonal polynomials on $[-1,1]$. As an application, the strong asymptotics of
these matrix orthogonal polynomials is derived, obtaining an explicit expression for
the corresponding Szeg\H o's matrix function.}

\vskip 0.25cm


\noindent {\it (1991) AMS Mathematics Subject Classification }: 42C05

\vskip 0.5cm

\noindent {\bf 1. Introduction. Semi-orthogonal functions}

\vskip 0.25cm

From long time ago, it is well known that there exists a simple relation between
orthogonal polynomials (OP) on the unit circle ($\T$) and OP on $[-1,1]$ (see [9, 24]).
This close relationship provides a method to translate results from OP on $\T$ to OP on
$[-1,1]$. For instance, this idea was largely exploited to get asymptotic properties of
OP on $[-1,1]$ starting from the asymptotics of OP on $\T$ [18, 19, 20, 24].
However, this
relation is valid only for symmetric measures on $\T$. Recently it has been shown that
above procedure can be generalized to arbitrary measures on $\T$, giving a connection
between any sequence of OP on $\T$ and the so-called semi-orthogonal functions [1, 4].

As we will see, semi-orthogonal functions are given in terms of a
sequence of two-dimensional matrix polynomials. The orthogonality
properties of semi-orthogonal functions implies that these matrix
polynomials are quasi-orthogonal with respect to some
two-dimensional matrix measure related to the measure on $\T$. A
sequence of matrix OP with respect to this matrix measure can be
explicitly constructed from above quasi-orthogonal matrix
polynomials. This gives a connection between OP on $\T$ and a
class of two-dimensional matrix OP on the real line.

Matrix OP on the real line appear in the Lanczos method  for block matrices [12, 13],
in the spectral theory of doubly infinite Jacobi matrices
[23] and discrete Sturm-Liouville operators [2, 3], in
the analysis of sequences of polynomials satisfying higher order recurrence relations
[8], rational approximation and system theory [10].
Unfortunately their study is much more complicated and few things are known if compared
with the scalar case (some nice surveys are [17, 21, 23]).

Previous connections between scalar and matrix OP appear in [14]
([15]) where it is derived a relation between scalar OP on an
algebraic harmonic curve (lemniscata) and matrix OP on the real
line (unit circle). Using a similar technique, a connection
between scalar OP with respect to a discrete Sobolev inner product
and matrix OP is presented in [8] (which is a consequence of the
fact that OP with respect to a discrete Sobolev inner product
satisfy a higher order recurrence relation, see [16]).  These
connections are more general than the one given in this paper
because they deal with matrix OP of arbitrary dimension. However,
they link matrix OP with ``unknown words", in the sense that not
too much is known about the different kind of OP that they connect
with matrix OP. So, they are not too useful to get new results for
matrix OP.

 On the contrary, the connection presented in this paper, although much more restricted,
let us translate results from the more ``known word" of scalar OP
on the unit circle to a great variety of two-dimensional matrix
OP. So, it provides many models of matrix OP where many things can
be known and that, therefore, can be used to get or check some
ideas about new results for matrix OP. Here we have to point out
that for certain applications, like the study of the doubly
infinite matrices that appear in discrete Sturm-Liouville problems
on the real line, only two-dimensional matrix OP are needed [3,
23].

 As an example of the utility of the present connection we derive the strong asymptotics
of these matrix OP when the corresponding matrix measure belongs
to the Szeg\H o's class. General results about this situation can
be found in [2], where a generalization of the connection between
the real line and $\T$ for matrix OP is used again to obtain the
asymptotics in the real line from the asymptotics in $\T$.
However, the problem is far from being closed since there is no
explicit expression for the Szeg\H o's matrix function that gives
the asymptotic behavior and only some general properties are
known. The connection given here let us obtain explicitly this
Szeg\H o's matrix function for a class of two-dimensional matrix
measures. Other results about asymptotics of matrix OP, such as
ratio and relative asymptotics, appear in [7] and [25]
respectively.

 Now, we proceed to introduce the starting point of our discussion, the semi-orthogonal
functions, summarizing some results in [1, 4] with a sketch of some proofs there
for the convenience of the reader.

  First of all we fix some notations. The real vector space of polynomials with
real coefficients is denoted by ${\cal P}$, the subspace of ${\cal
P}$ of polynomials with degree less than or equal to $n$ is ${\cal
P}_n$ and ${\cal P}_n^\#$ is the subset of ${\cal P}_n$
constituted by those polynomials whose degree is exactly $n$.
Also, $\Lambda$ is  the complex vector space of Laurent
polynomials, that is, $\Lambda = \bigcup_{n=0}^\infty
\Lambda_{-n,n}$ where $\Lambda_{m,n} = \left\{ \sum_{k=m}^n
\alpha_k z^k \big| \alpha_k \in \C \right\}$ for $m \leq n$. The
elements of $\Lambda_{m,n}$ such that $\alpha_m, \alpha_n \neq 0$
form the subset $\Lambda_{m,n}^\#$. For an arbitrary complex
number $\alpha$, their real and imaginary parts are denoted $\Re
\alpha$ and $\Im \alpha$ respectively.

 Taking into account the usual identification between the unit circle
$\T=\{e^{i\theta}| \; \theta \in [0,2\pi)\}$ and the interval
$[0,2\pi)$, we talk about a measure on $\T$ when we deal with a
measure with support on $[0,2\pi)$. With this convention, in what
follows $d\mu$ is a measure on $\T$ with finite moments. Unless we
say explicitly that it is an arbitrary measure on $\T$, we suppose
that $d\mu$ is a positive measure with infinite support. Then, the
sesquilinear functional $\big<\cdot,\cdot\big>_{d\mu}$ on
$\Lambda$ defined by
$$
\big< f,g \big>_{d\mu} = \int_0^{2\pi}
f(e^{i\theta}) \overline{g(e^{i\theta})} \, d\mu(\theta),
\quad f,g \in \Lambda,
$$
is an inner product and, hence, there exists a unique sequence $\sn \phi 0$ of
monic OP with respect to $\big<\cdot,\cdot\big>_{d\mu}$. If, as it is usual, $\phi_n^*$
denotes the  reversed polynomial of
$\phi_n$ ($\phi_n^*(z) = z^n \overline{\phi}_n(z^{-1})$), then, it is well known that
OP are determined by the so-called Schur parameters $a_n = \phi_n(0)$ through the
recurrence
$$
\eqalign{
& \phi_0(z) = 1, \cr
& \phi_n(z) = z \phi_{n-1}(z) + a_n \phi_{n-1}^*(z), \quad n \geq 1. \cr
}
\eqno (1)
$$
 If we denote by $b_n$ the coefficient of $z^{n-1}$ in $\phi_n(z)$, from (1) we
have that
$$
b_n = b_{n-1} + a_n \overline{a}_{n-1}, \quad n \geq 1.
\eqno (2)
$$
 Notice that $b_0=0$ and
$$
b_n = \sum_{k=1}^{n} a_k \overline{a}_{k-1}, \quad n \geq 1.
\eqno (3)
$$
 We can use (1) to show that the positive constants
$\varepsilon_n = \big< \phi_n,\phi_n \big>_{d\mu}$ are related to the Schur parameters
by
$$
{\varepsilon_n \over \varepsilon_{n-1}} = 1-|a_n|^2, \quad n \geq 1.
\eqno (4)
$$

 This relation implies that $|a_n|<1$ for $n \geq 1$ and that the sequence
$\sn \varepsilon 0$ must be strictly decreasing. Besides, (4)
gives the following expression for $\varepsilon_n$
$$
\varepsilon_n = \prod_{k=1}^n (1-|a_k|^2) \varepsilon_0, \quad n \geq 1.
\eqno (5)
$$

  Orthonormal polynomials are defined up to a factor with unit module, but they can be
fixed if we ask for their leading coefficients to be real and
positive. In this case we denote the $n$-th orthonormal polynomial
by $\varphi_n$, and the corresponding leading coefficient by
$\kappa_n$. It is clear that $\kappa_n = \varepsilon_n^{-1/2}$
and, thus, $\sn \kappa 0$ is strictly increasing.

  The symmetric measure of $d\mu$ is
$$
d\widetilde\mu(\theta) = -d\mu(2\pi - \theta), \quad \theta \in [0,2\pi),
$$
and the measure $d\mu$ is said to be symmetric iff $d\widetilde\mu=d\mu$. This is
equivalent to affirm that the monic OP have real coefficients, which, in sight of (1),
is in fact equivalent to state that the Schur parameters are real.

  With the intention of connecting $\T$ with the interval $[-1,1]$, for
$z \in \C \setminus \{0\}$ we write $x=(z+z^{-1})/2$ and
$y=(z-z^{-1})/2i$ (therefore $z=x+iy, z^{-1}=x-iy$ and
$x^2+y^2=1$). Both expressions give a transformation in the
complex plane that maps $\T$ on the interval $[-1,1]$. Moreover,
they map bijectively onto $\C \setminus [-1,1]$ the exterior of
$\T$ as well as its interior excepting the origin. So, when
restricted to these domains we can invert the transformations
giving, for example, $z = x + \sqrt{x^2-1}$ (the choice of the
square root must be done according to the location of $z$:
exterior or interior to $\T$). Also, the transformation
$x=(z+z^{-1})/2$ maps biyectively the upper as well as the lower
closed half $\T$ onto $[-1,1]$ (in this case, writing
$z=e^{i\theta}$, it is $x=\cos\theta$). So, by composition with
the corresponding inverse transformations, the measure $d\mu$
provides two projected measures $d\nu_1$, $d\nu_2$ on $[-1,1]$,
being
$$
\eqalign{
& d\nu_1(x) = -d\mu(\arccos x), \cr
& d\nu_2(x) = -d\widetilde\mu(\arccos x). \cr
}
\eqno (6)
$$
The condition of symmetry for $d\mu$ is equivalent to the equality $d\nu_1 = d\nu_2$.

    Now, we wish to arrive at a family of polynomials with real coefficients, orthogonal
with respect to an inner product defined through the measures $d\nu_1$, $d\nu_2$. To
this end, and following [1, 4], we start by introducing previously the so called
semi-orthogonal functions.

\noindent {\bf Definition 1.}
The semi-orthogonal functions (SOF) associated to the measure $d\mu$ are the functions
${f^{(k)}_n} \colon \C \setminus \{0\} \to \C, \; n \geq 1, \; k=1,2$, defined by
$$
\eqalign{
& f^{(1)}_n (z) = {z \phi_{2n-1}(z) + \phi_{2n-1}^*(z) \over 2^n z^n}, \cr
& f^{(2)}_n (z) = {z \phi_{2n-1}(z) - \phi_{2n-1}^*(z) \over i 2^n z^n}, \cr
}
$$
where $\phi_n, \; n \geq 1$, are the monic OP with respect to
$\big<\cdot,\cdot\big>_{d\mu}$.

 The expressions in above definition are the same used in Szeg\H o's method, with the
difference that we consider here monic OP with complex instead of real coefficients.
Let us go to summarize some interesting properties of SOF [1, 4]:

\noindent {\bf Proposition 1.}
{\sl The SOF associated to $d\mu$ satisfy:

\item{i)} $\overline{f}^{(k)}_n (z^{-1}) = f^{(k)}_n (z)$ and there is a unique
decomposition
$$
f^{(k)}_n (z) = f^{(k1)}_n (x) + y \, f^{(k2)}_n (x), \quad f^{(kj)}_n \in {\cal P}.
$$
More precisely,
$f^{(11)}_n \in {\cal P}_n^\#$, $f^{(22)}_n \in {\cal P}_{n-1}^\#$, both monic
polynomials, and
$f^{(21)}_n \in {\cal P}_{n-1}, f^{(12)}_n \in {\cal P}_{n-2}$.

\item{ii)} The family of functions
${\cal B}_n = \{1\} \cup \left( \bigcup_{m=1}^n \{f^{(1)}_m, f^{(2)}_m\} \right)$
is a basis of $\Lambda_{-n,n}$ for all $n \geq 1$. The matrix of
$\big<\cdot,\cdot\big>_{d\mu}$ with respect to the basis
${\cal B} = \cup_{n \geq 1} {\cal B}_n$ of $\Lambda$ is a diagonal-block one
$$
\pmatrix{
\varepsilon_0 & 0 & 0 & \dots \cr
0 & C_1 & 0 & \dots \cr
0 & 0 & C_2 & \dots \cr
\vdots & \vdots & \vdots & \ddots \cr
},
\eqno (7)
$$
where
$$
C_n = {\varepsilon_{2n-1} \over 2^{2n-1}}
\pmatrix{1-\Re a_{2n} & -\Im a_{2n} \cr -\Im a_{2n} & 1+\Re a_{2n} \cr},
\quad n \geq 1,
\eqno (8)
$$
\item{} being $a_n$ the Schur parameters related to $d\mu$ and $\varepsilon_n$ given in
(5).
}

\noindent {\bf Proof.} From the definition of $\phi_n^*$ we have that
$$
\eqalign{
& f^{(1)}_n(z) = 2^{-n}
({z^{-n+1} \phi_{2n-1}(z)} + z^{n-1} \overline{\phi}_{2n-1}(z^{-1})), \cr
& f^{(2)}_n(z) = -i2^{-n}
({z^{-n+1} \phi_{2n-1}(z)} - z^{n-1} \overline{\phi}_{2n-1}(z^{-1})), \cr
}
$$
and, thus, $\overline{f}^{(k)}_n (z^{-1}) = f^{(k)}_n (z)$.

 If the decomposition given in i) exits, it must be unique. If we suppose two such
decompositions
$$
f^{(k)}_n(z) = f^{(k1)}_n (x) + y \, f^{(k2)}_n (x)
= g^{(k1)}_n (x) + y \, g^{(k2)}_n (x),
\quad f^{(kj)}_n, g^{(kj)}_n \in {\cal P},
$$
then
$$
f^{(k)}_n(z^{-1}) = f^{(k1)}_n (x) - y \, f^{(k2)}_n (x)
= g^{(k1)}_n (x) - y \, g^{(k2)}_n (x),
$$
and above equalities give $f^{(kj)}_n = g^{(kj)}_n$ for $k=1,2$.

 To see that this decomposition exist, let us write
$\phi_n(z) = \sum_{k=0}^n \alpha_k z^k, \; \alpha_k \in \C$, with $\alpha_n=1$. Then,
$$
\eqalign{
f^{(1)}_n (z) & = 2^{-n} z^{-n} \sum_{k=0}^{2n-1}
                 (\alpha_k z^{k+1} + \overline{\alpha}_k z^{2n-k-1}) \cr
& = T_n(x) + \sum_{j=0}^{n-1} 2^{j-n} \Re (\alpha_{n-j-1}+\alpha_{n+j-1}) T_j(x) \cr
& \quad + y \sum_{j=0}^{n-2} 2^{j+1-n} \Im (\alpha_{n-j-2}+\alpha_{n+j}) U_j(x), \cr
}
\eqno (9)
$$
where $T_j(x)=2^{-j}(z^j+z^{-j})$ and
$U_j(x)= -i y^{-1} 2^{-j-1}(z^{j+1}-z^{-j-1})$ are respectively the $j$-th Tchebychev
monic polynomials of first and second kind (see for example [5] or [24]). This proves
i) for $f^{(1)}_n$. The proof for $f^{(2)}_n$ is similar.

 Notice that $f^{(k)}_n \in \Lambda_{-n,n}^\#$ for $i=1,2$. In fact,
$f^{(1)}_n(z) = 2^{-n} (z^n + z^{-n}) + \dots$ and
$f^{(2)}_n(z) = -i2^{-n} (z^n - z^{-n}) + \dots$, where the dots mean terms belonging to
$\Lambda_{-n+1,n-1}$. Thus, it is obvious that ${\cal B}_n$ is a basis of
$\Lambda_{-n,n}$. The block-diagonal structure (7) of the matrix of
$\big<\cdot,\cdot\big>_{d\mu}$ with respect to ${\cal B},$
$$
\pmatrix{
\big<1,1\big>_{d\mu} & \big<1,f^{(1)}_1\big>_{d\mu} & \big<1,f^{(2)}_1\big>_{d\mu} &
\dots \cr
\big<f^{(1)}_1,1\big>_{d\mu} & \big<f^{(1)}_1,f^{(1)}_1\big>_{d\mu} &
\big<f^{(1)}_1,f^{(2)}_1\big>_{d\mu} & \dots \cr
\big<f^{(2)}_1,1\big>_{d\mu} & \big<f^{(2)}_1,f^{(1)}_1\big>_{d\mu} &
\big<f^{(2)}_1,f^{(2)}_1\big>_{d\mu} & \dots \cr
\vdots & \vdots & \vdots & \ddots \cr
}
$$
is just a direct consequence of the orthogonality relations for $\phi_n$ and
$\phi_n^*$, that is, the only conditions
$\big< \phi_n,z^k \big>_{d\mu} = \big< \phi_n^*,z^{k+1} \big>_{d\mu} = 0$ for
$0 \leq k \leq n-1$ imply
$\big<1,f^{(k)}_n\big>_{d\mu} = \big<f^{(k)}_n,f^{(j)}_m\big>_{d\mu} = 0$ for
$n \neq m$ and $k,j=1,2$.

Finally, the expression (8) for the matrix
$$
C_n = \pmatrix{
\big<f^{(1)}_n,f^{(1)}_n\big>_{d\mu} & \big<f^{(1)}_n,f^{(2)}_n\big>_{d\mu} \cr
\big<f^{(2)}_n,f^{(1)}_n\big>_{d\mu} & \big<f^{(2)}_n,f^{(2)}_n\big>_{d\mu} \cr
}
$$
follows straightforward from the relations $\big<
\phi_n^*,\phi_n^* \big>_{d\mu} = \big< \phi_n,\phi_n \big>_{d\mu}
= \varepsilon_n$ and $\big< z\phi_n,\phi_n^* \big>_{d\mu} = -
\varepsilon_n a_{n+1}$, the last one obtained from the recurrence
formula (1). $\qed$

 Notice that the property i) implies that SOF are real on $\T$.
The incomplete orthogonality expressed in ii) of previous proposition is the origin
of the name \lq\lq semi-orthogonal functions" given to the functions $f^{(k)}_n$. Notice
that, when the measure $d\mu$ is symmetric, the monic OP $\phi_n$ have real coefficients
and it follows from previous proof that
$f^{(12)}_n = f^{(21)}_n = 0, \; n \geq 1$. Moreover, in this case the Schur
parameters are real and, then, the SOF are indeed strictly orthogonal.

Before continuing, it is useful to introduce a new notation.

\noindent{\bf Definition 2.}
The vector semi-orthogonal functions (VSOF) associated to the measure $d\mu$ are the
functions
$\b f_n \colon \C \setminus \{0\} \to \C^2, \; n \geq 0$, defined by
$$
\b f_n(z) = \pmatrix{f^{(1)}_n(z) \cr f^{(2)}_n(z) \cr},
$$
where $f^{(1)}_0(z)=1, f^{(2)}_0(z)=0$ and $f^{(k)}_n(z), \; n \geq 1, \; k=1,2$, are
the SOF related to $d\mu$.

\noindent {\bf Remark 1.}
Proposition 1 i), that is given for $n \geq 1$, holds for $n=0$ too, being
$f^{(11)}_0(x)=1,f^{(12)}_0(x)=f^{(21)}_0(x)=f^{(22)}_0(x)=0$. We will refer to
$(f^{(k)}_n)_{n \geq 0 \atop k=1,2}$, as the complete family of SOF associated to
$d\mu$.

\noindent {\bf Remark 2.}
 The fact that ${\cal B}_n$ is a basis of $\Lambda_{-n,n}$ ensures that
$({\b f_m})_{m=0}^n$ is a set of generators for the modulus $\Lambda_{-n,n}^2$ over
the ring $\C^{(2,2)}$ of $2 \times 2$ complex matrices. Hence, $({\b f_m})_{m \geq
0}$ is a set of generators for $\Lambda^2$.
Although $({\b f_n})_{n \geq 0}$ is not a basis, it is not difficult to
see that in the decomposition of an arbitrary element of $\Lambda^2$ as a linear
combination of $({\b f_n})_{n \geq 0}$, all the matrix coefficients are
univocally determined excepting the one related to $\b f_0$, whose first column is
determined whereas the second one is arbitrary.

\noindent{\bf Definition 3.}
Given an arbitrary measure $d\mu$ on $\T$ we define the sesquilinear functional
$\;\;\lp\cdot,\cdot\rp_{d\mu} \colon \Lambda^2 \times \Lambda^2 \to \C^{(2,2)}$
in the following way
$$
\lp \b f,\b g \rp_{d\mu} =
\int_0^{2\pi} \b f(e^{i\theta}) \b g(e^{i\theta})^* \, d\mu(\theta), \quad
\b f, \b g \in \Lambda^2,
$$
where, for an arbitrary matrix $A$, we write $A^* = \overline{A}^T$ and the symbol $T$
denotes the operation of transposition.

\noindent {\bf Remark 3.}
If $\b f = \pmatrix{f^{(1)} \cr f^{(2)} \cr},
\b g = \pmatrix{g^{(1)} \cr g^{(2)} \cr}$ with $f^{(k)},g^{(k)} \in \Lambda$ for
$k=1,2$, then
$$
\lp \b f,\b g \rp_{d\mu} =
\pmatrix{
\big<f^{(1)},g^{(1)}\big>_{d\mu} & \big<f^{(1)},g^{(2)}\big>_{d\mu} \cr
\big<f^{(2)},g^{(1)}\big>_{d\mu} & \big<f^{(2)},g^{(2)}\big>_{d\mu} \cr
}.
$$
Notice that for all $\b f,\b g \in \Lambda^2$
it is $\lp z\b f,z\b g\rp_{d\mu} = \lp\b f,\b g\rp_{d\mu}$ and
$\lp\b g,\b f\rp_{d\mu} = \lp\b f,\b g\rp_{d\mu}^*$.

The orthogonality properties of SOF, translated to the language of VSOF, give the
following result.

\noindent {\bf Proposition 2.}
{\sl The VSOF associated to $d\mu$ are orthogonal with respect to
$\lp\cdot,\cdot\rp_{d\mu}$. More precisely,
$$
\lp \b f_n, \b f_m \rp_{d\mu} = C_n \delta_{n,m}, \quad n,m \geq 0,
$$
where $C_n$ is given in (8) for $n \geq 1$ and $C_0 = \varepsilon_0 C$ with
$C = \pmatrix{1 & 0 \cr 0 & 0 \cr}$.
}

\noindent{\bf Definition 4.}
The Schur matrices associated to the measure $d\mu$ are the following real symmetric
traceless matrices
$$
H_n = \pmatrix{\Re a_n & \Im a_n \cr \Im a_n & -\Re a_n \cr}, \quad n \geq 0,
$$
where $a_n$ are the Schur parameters related to $d\mu$.

\noindent {\bf Remark 4.}
Notice that we can write
$$
C_n = {\varepsilon_{2n-1} \over 2^{2n-1}} (I-H_{2n}), \quad n \geq 1.
\eqno (10)
$$
For $n \geq 1$, the condition $|a_n|<1$, which is equivalent to $|\det H_n| < 1$,
ensures that $C_n$ is positive definite and, therefore, nonsingular. With above
notation,
$$
C_n^{-1} = {2^{2n-1} \over \varepsilon_{2n}} (I+H_{2n}), \quad n \geq 1,
\eqno (11)
$$
where we have used (4).

Bearing in mind that $(\b f_m)_{m=0}^n$ is a set of generators for
$\Lambda_{-n,n}$, we get from Proposition 2 the following consequence.

\noindent {\bf Corollary 1.}
{\sl The VSOF associated to $d\mu$ satisfy for $n \geq 1$
$$
\lp \b f_n, \b f \rp_{d\mu} = 0, \quad \forall \b f \in \Lambda_{-n+1,n-1}^2.
$$
}

\vskip 0.3cm

\noindent {\bf 2. Recurrence relation for semi-orthogonal functions}

\vskip 0.25cm

The VSOF associated to a measure form a set of orthogonal vector Laurent polynomials,
where the orthogonality is respect to some sesquilinear functional related to the measure.
The natural question that arises is if, analogously to OP, they satisfy a three-term
recurrence relation. The answer to this question is given in the following proposition.

\noindent {\bf Proposition 3.}
{\sl The VSOF associated to $d\mu$ satisfy the recurrence relation
$$
z \b f_n(z) =  (I+iJ) \b f_{n+1}(z) + \b L_n \b f_n(z) + \b M_n \b f_{n-1}(z),
\quad n \geq 1,
$$
where $I$ is the $2\times2$ identity matrix and
$$
\eqalign{
& J=\pmatrix{0 & 1 \cr -1 & 0}, \cr
& \b L_n = {1 \over 2}
\left\{ (I - H_{2n}) H_{2n-1} (I + iJ) - (I + iJ) H_{2n+1} (I + H_{2n}) \right\},
\quad n \geq 1, \cr
& \b M_n = {1 \over 4} \det (I-H_{2n-1}) (I - H_{2n}) (I - iJ) (I + H_{2n-2}),
\quad n \geq 1, \cr
}
$$
being $H_n$ the Schur matrices related to $d\mu$.
}

\noindent {\bf Proof.} As it is usual, we begin by decomposing $z \b f_n$ with
respect to the set of generators $({\b f_n})_{n \geq 0}$. Since
$z\b f_n \in \Lambda_{-n+1,n+1}^2 \subset \Lambda_{-n-1,n+1}^2$, it is obvious
that
$$
z \b f_n(z) = \sum_{k=0}^{n+1} A^{(n)}_k \b f_k(z),
\quad A^{(n)}_k \in \C^{(2,2)}, \quad n \geq 1.
$$
From Remark 2 we know that the matrix coefficients $A^{(n)}_k$ are univocally
determined for $k \geq 1$, while for $A^{(n)}_0$ only the product $A^{(n)}_0 \b f_0$ is
fixed, that is, the first column of $A^{(n)}_0$ is determined whereas the second one is
arbitrary.

Now, by projecting above decomposition of $z \b f_n$ over $\b
f_j$, with $j \leq n+1$, we find that $\lp z\b f_n , \b f_j
\rp_{d\mu} = A^{(n)}_j C_j$.

When $j=0$ and $n \geq 2$ we get
$A^{(n)}_0 C_0 = \lp \b f_n , z^{-1}\b f_0 \rp_{d\mu} = 0$ due to Corollary 1.
Therefore,
$A^{(n)}_0 \b f_0 = A^{(n)}_0 C \b f_0 = 0$.

If $1 \leq j \leq n-2$, then $A^{(n)}_j C_j = \lp \b f_n,\;\;
z^{-1}\b f_j \rp_{d\mu} = 0$ again by means of Corollary 1. Now,
the regularity of $C_j$ for $j \geq 1$ forces $A^{(n)}_j = 0$.

Hence, we can write
$$
z \b f_n(z) =
A^{(n)}_{n+1} \b f_{n+1}(z) + A^{(n)}_n \b f_n(z) + A^{(n)}_{n-1} \b f_{n-1}(z),
\quad n \geq 1,
\eqno (12)
$$
where the only indetermination is in the second column of
$A^{(1)}_0$ that can be arbitrarily chosen. Notice that the
coefficients $A^{(n)}_{n-1}$ and $A^{(n-1)}_n$ are related for $n
\geq 1$ by $A^{(n)}_{n-1} C_{n-1} = \lp z\b f_n , \b f_{n-1}
\rp_{d\mu} = \lp z^{-1} \b f_{n-1} , \b f_n \rp_{d\mu}^* = \lp z\b
f_{n-1} , \b f_n \rp_{d\mu}^T = C_n (A^{(n-1)}_n)^T$, where we
have used the fact that VSOF, like SOF, are real on $\T$. Thus,
$$
\eqalign{
& A^{(1)}_0 \b f_0 = \varepsilon_0^{-1} A^{(1)}_0 C_0 \b f_0
= \varepsilon_0^{-1} C_1 (A^{(0)}_1)^T \b f_0 \cr &
A^{(n)}_{n-1} = C_n (A^{(n-1)}_n)^T C_{n-1}^{-1} \cr
}
\eqno (13)
$$
This means that we only need to calculate $A^{(n)}_{n+1}$ and $A^{(n)}_{n}$.
To this end we introduce the following elements of $\Lambda^2$
$$
\eqalign{
& \b g_0(z) = \b f_0(z) \cr
& \b g_n(z) = \pmatrix{z^n \cr z^{-n}}, \quad n \geq 1, \cr
}
$$
so that, $(\b g_n)_{n \geq 0}$ is a set of generators for the module $\Lambda^2$ with
decomposition properties similar to those above described for
$(\b f_n)_{n \geq 0}$.
Thus, we can decompose both sides of (12) in
$(\b g_n)_{n \geq 0}$ and, then, equal coefficients of $\b g_n$ for $n \geq 1$.

  We begin with the decomposition of $\b f_n$ for $n \geq 1$.
If $Q = \pmatrix{1 & 1 \cr -i & i \cr}$, then
$$
\eqalign{
\b f_n(z) &= 2^{-n} Q
\pmatrix{z^{1-n} \phi_{2n-1}(z) \cr z^{-n} \phi_{2n-1}^*(z) \cr} \cr
&= 2^{-n} Q \left\{ \b g_n(z) +
\pmatrix{b_{2n-1} & a_{2n-1} \cr \overline{a}_{2n-1} & \overline{b}_{2n-1} \cr}
\b g_{n-1}(z) + \dots \right\}, \cr
}
\eqno (14)
$$
where the dots mean terms belonging to $\Lambda_{-n+2,n-2}^2$ for
$n \geq 2$ and no terms for $n=1$. Besides, we need the
decomposition of $z \b f_n$ for $n \geq 1$
$$
\eqalign{
z \b f_n(z) &= 2^{-n} Q
\pmatrix{z^{2-n} \phi_{2n-1}(z) \cr z^{1-n} \phi_{2n-1}^*(z) \cr} \cr
&= 2^{-n} Q \left\{ C_0 \b g_{n+1}(z) +
\pmatrix{b_{2n-1} & a_{2n-1} \cr \overline{a}_{2n-1} & \overline{b}_{2n-1} \cr}
C_0 \b g_n(z) + \dots \right\}, \cr
}
\eqno (15)
$$
where now the dots mean terms belonging to $\Lambda_{-n+1,n-1}^2$.

Introducing (14) and (15) into (12) and equaling coefficients of $\b g_{n+1}$ and
$\b g_n,$ lead to
$$
\eqalign{
Q C_0 & = {1 \over 2} A^{(n)}_{n+1} Q, \cr
Q \pmatrix{b_{2n-1} & a_{2n-1} \cr \overline{a}_{2n-1} & \overline{b}_{2n-1} \cr} C_0
& = {1 \over 2} A^{(n)}_{n+1} Q
\pmatrix{b_{2n+1} & a_{2n+1} \cr \overline{a}_{2n+1} & \overline{b}_{2n+1} \cr}
+ A^{(n)}_n Q,
}
$$
which have the solutions
$$
\eqalign{
A^{(n)}_{n+1} & = 2 Q C_0 Q^{-1}, \cr
A^{(n)}_n & = Q
\pmatrix{b_{2n-1} - b_{2n+1} & -a_{2n+1} \cr \overline{a}_{2n-1} & 0 \cr} Q^{-1} \cr
& = Q \left\{
\overline{a}_{2n-1} \pmatrix{-a_{2n} & 0 \cr 1 & 0 \cr} -
a_{2n+1} \pmatrix{\overline{a}_{2n} & 1 \cr 0 & 0 \cr}
\right\} Q^{-1} \cr
& = Q \left\{
\overline{a}_{2n-1}
\pmatrix{1 & -a_{2n} \cr -\overline{a}_{2n} & 1 \cr}
\pmatrix{0 & 0 \cr 1 & 0 \cr} -
a_{2n+1}
\pmatrix{0 & 1 \cr 0 & 0 \cr}
\pmatrix{1 & a_{2n} \cr \overline{a}_{2n} & 1 \cr}
\right\} Q^{-1}. \cr
}
$$
Here we have used relation (2) between $a_n$ and $b_n$. At this
point it is useful to notice that
$$
\eqalign{
& Q C_0 Q^{-1} = {1 \over 2} (I + iJ), \quad J = \pmatrix{0 & 1 \cr -1 & 0 \cr} \cr
& Q \pmatrix{1 & \pm a_n \cr \pm \overline{a}_n & 1 \cr} Q^{-1} = I \pm H_n,
\quad n \geq 0, \cr
& \overline{a}_n Q \pmatrix{0 & 0 \cr 1 & 0 \cr} Q^{-1} = {1 \over 2} H_n (I + iJ),
\quad n \geq 0, \cr
& a_n Q \pmatrix{0 & 1 \cr 0 & 0 \cr} Q^{-1} = {1 \over 2} (I + iJ) H_n,
\quad n \geq 0, \cr
}
$$
where $H_n$ are the Schur matrices related to $d\mu$.
Thus, we finally find
$$
\eqalign{
& A^{(n)}_{n+1} = I + iJ, \quad n \geq 1, \cr
& A^{(n)}_n = {1 \over 2}
\left\{ (I - H_{2n}) H_{2n-1} (I + iJ) - (I + iJ) H_{2n+1} (I + H_{2n}) \right\},
\quad n \geq 1. \cr
}
\eqno (16)
$$

    Now, by using (4), (10), (11), (13) and (16), we see that
$$
\eqalign{
& A^{(1)}_0 \b f_0 = {1 \over 2} (1 - |a_1|^2) (I - H_2) (I - iJ) \b f_0, \cr
& A^{(n)}_{n-1} = {1 \over 4} (1 - |a_{2n-1}|^2) (I - H_{2n}) (I - iJ) (I + H_{2n-2}),
\quad n \geq 2. \cr
}
$$
It is possible to choose arbitrarily the second column of $A^{(1)}_0$,
so we can fixe it to be null, that is
$$
A^{(1)}_0 = {1 \over 2} (1 - |a_1|^2) (I - H_2) (I - iJ) C =
{1 \over 4} (1 - |a_1|^2) (I - H_2) (I - iJ) (I + H_0).
$$
Hence, all the coefficients $A^{(n)}_{n-1}$ can be given by
$$
A^{(n)}_{n-1} = {1 \over 4} \det (I-H_{2n-1}) (I - H_{2n}) (I - iJ) (I + H_{2n-2}),
\quad n \geq 1.
$$
Since $\b L_n=A^{(n)}_n$ and $\b M_n=A^{(n)}_{n-1}$ the proposition
is proved.
$\qed$

\vskip 0.5cm

\noindent {\bf 3. Semi-orthogonal functions and matrix measures}

\vskip 0.25cm

Now we are going to translate previous results on $\Lambda^2$ to the modulus
${\cal P}^{(2,2)}$ of $2 \times 2$ matrices with coefficients in ${\cal P}$.
To do that we will associate a sequence of matrix polynomials to any family of SOF.

\noindent {\bf Definition 5.}
The matrix polynomials $\b F_n \in {\cal P}^{(2,2)}, \; n \geq 0$, associated to
the measure $d\mu,$ are defined by
$$
\b F_n(x) =
\pmatrix{f^{(11)}_n(x) & f^{(12)}_n(x) \cr f^{(21)}_n(x) & f^{(22)}_n(x) \cr},
$$
where
$f^{(k)}_n(z)=f^{(k1)}_n(x)+yf^{(k2)}_n(x), \, n \geq 0, \; k=1,2$, is the
unique decomposition given in Proposition 1 and Remark 1 for the complete family of SOF
related to $d\mu$.

\noindent {\bf Remark 5.}
The VSOF and the matrix polynomials associated to $d\mu$ are related by
$$
\b f_n(z) = \b F_n(x) \pmatrix{1 \cr y}, \quad n \geq 0,
\eqno (17)
$$
and, thus,
$$
\pmatrix{\b f_n(z) & \b f_n(z^{-1})} = \b F_n(x) \pmatrix{1 & 1 \cr y & -y},
\quad n \geq 0,
$$
from which we get
$$
\b F_n(x) = {1 \over 2} \pmatrix{\b f_n(z) & \b f_n(z^{-1})}
\pmatrix{1 & y^{-1} \cr 1 & -y^{-1}}, \quad n \geq 0,
\eqno (18)
$$
for $x \neq \pm 1$.

It is natural to expect for the above matrix polynomials to inherit some orthogonality
properties from the ones satisfied by the corresponding SOF. To see this we introduce
the following matrix measure.

\noindent {\bf Definition 6.}
Given an arbitrary measure $d\mu$ on $\T$ the matrix measure $d\Omega$ associated to
$d\mu$ is the following $2\times2$ symmetric matrix measure on $[-1,1]$
$$
d\Omega(x) = {1 \over 2}
\pmatrix{d\rho(x) & \sqrt{1-x^2} d\sigma(x) \cr
\sqrt{1-x^2} d\sigma(x) & (1-x^2) d\rho(x) \cr},
\eqno (19)
$$
where $d\rho$ and $d\sigma$ are scalar measures on $[-1,1]$ given by
$$
\eqalign{
& d\rho(x) = d\nu_1(x) + d\nu_2(x), \cr
& d\sigma(x) = d\nu_1(x) - d\nu_2(x), \cr
}
$$
and $d\nu_1$, $d\nu_2$ are the projected measures of $d\mu$ defined in (6).

\noindent {\bf Remark 6.}
Notice that a matrix measure $d\Omega$ with the form (19), being $d\rho$ and $d\sigma$
arbitrary scalar measures on $[-1,1]$, is always associated to some measure $d\mu$ on
$\T$. The related measure $d\mu$ is positive iff $d\nu_1$, $d\nu_2$ so are, which
holds iff $|d\sigma| \leq d\rho$ (this implies that $d\rho$ is positive and that
supp$(d\sigma) \subset$ supp$(d\rho)$). Therefore, when $d\mu$ is positive it
has an infinite support iff $d\rho$ so does.

Now, the results in Proposition 1 for SOF have the following consequences for the
corresponding matrix polynomials.

\noindent {\bf Proposition 4.}
{\sl The matrix polynomials $(\b F_n)_{n\geq 0}$ and the matrix measure $d\Omega$
associated to $d\mu$ satisfy:

\item{i)} $\deg \b F_n = n$. More precisely, $\b F_0(x)=C$, with the matrix $C$ as in
Proposition 2, and
$$
\b F_{n+1}(x) = C x^{n+1} + \pmatrix{\eta_n & 0 \cr \gamma_n & 1} x^n + \dots,
\quad n \geq 0,
$$
where the dots mean terms with degree less than $n$ and
$$
\eta_n = {1 \over 2} \Re (a_{2n+1}+b_{2n+1}), \quad
\gamma_n = {1 \over 2} \Im (a_{2n+1}+b_{2n+1}),
$$
being $a_n$ the Schur parameters related to $d\mu$ and $b_n$ given in (3),

\item{ii)} $\int_{-1}^1 \b F_n(x) \, d\Omega(x) \, \b F_m^T(x) =
{1 \over 2} C_n \delta_{n,m}, \quad n,m \geq 0$ .
}

\noindent {\bf Proof.}
The result i) follows straightforward from Proposition 1 i) and Definition 5,
wherefrom we see that $\gamma_n$ is the leading coefficient of $f^{(21)}_{n+1}(x)$,
while $\eta_n$ is the coefficient of $x^n$ in $f^{(11)}_{n+1}(x)$. The expression (9)
for $f^{(1)}_n$ shows that $\eta_n = {1 \over 2} \Re (a_{2n+1}+b_{2n+1})$ and a similar
expression for $f^{(2)}_n$ gives $\gamma_n = {1 \over 2} \Im (a_{2n+1}+b_{2n+1})$.

 To prove ii) it is enough to notice that, using (17), we get from
Definition 3
$$
\eqalign{
\lp \b f_n,\b f_m \rp_{d\mu}
& = \int_0^{2\pi} \b f_n(e^{i\theta}) \b f_m^T(e^{i\theta})
\, d\mu(\theta) \cr
& =  \int_{-1}^1 \b F_n(x)
\pmatrix{1 & \sqrt{1-x^2} \cr \sqrt{1-x^2} & 1-x^2}
\b F_m^T(x) \, d\nu_1(x) \cr
& \quad + \int_{-1}^1 \b F_n(x)
\pmatrix{1 & -\sqrt{1-x^2} \cr -\sqrt{1-x^2} & 1-x^2}
\b F_m^T(x) \, d\nu_2(x), \cr
}
$$
whit the positive choice for the square root. Taking into account
Definition 6 we see that
$$
\lp \b f_n,\b f_m \rp_{d\mu} = 2\int_{-1}^1 \b F_n(x) \, d\Omega(x) \, \b F_m^T(x),
$$
and Proposition 2 gives ii).
$\qed$

\noindent {\bf Remark 7.}
From Proposition 4 i) we see that
$$
(I-C) \b F_{n+1}(x) + \b F_n(x) = \Gamma_n x^n + \dots,
\quad \Gamma_n = \pmatrix{1 & 0 \cr \gamma_n & 1 \cr}, \quad n \geq 0,
\eqno (20)
$$
where, again, the dots mean terms with degree less than $n$. So, it is
obvious that every element of ${\cal P}_n^{(2,2)}$ is a linear combination of
$(\b F_m)_{m=0}^{n+1}$, and, therefore, $(\b F_n)_{n \geq 0}$ is a set of
generators for ${\cal P}^{(2,2)}$.

 Unfortunately, in spite of Proposition 4 ii), we can not say that $(\b F_n)_{n\geq 0}$
is a sequence of left orthogonal matrix polynomials with respect to $d\Omega$. A
sequence $(\b P_n)_{n\geq 0}$ of $2\times2$ real matrix polynomials is called a sequence
of left orthogonal matrix polynomials (LOMP) with respect to $d\Omega$ if it satisfies
[6, 17, 23]:

\item{(I)} $\deg \b P_n = n$, and the leading coefficient of $\b P_n$ is nonsingular.

\item{(II)} $\int_{-1}^1 \b P_n(x) \, d\Omega(x) \, x^k = 0$ for $0 \leq k < n$ and
$\int_{-1}^1 \b P_n(x) \, d\Omega(x) \, x^n$ is nonsingular.

\noindent However, the leading coefficient $C$ of $\b F_n$ is singular. Even more,
although
$$
\int_{-1}^1 \b F_n(x) \, d\Omega(x) \, x^k = 0, \quad 0 \leq k \leq n-2,
\eqno (21)
$$
since $\b F_n$ is orthogonal to
span$\{ \b F_0, \b F_1, \dots , \b F_{n-1} \} \supset {\cal P}_{n-2}^{(2,2)}$,
we have from (20) that
$$
\int_{-1}^1 \b F_n(x) \, d\Omega(x) \, x^{n-1} =
{1 \over 2} C_n  (I-C)^T (\Gamma_{n-1}^{-1})^T = {1 \over 2} C_n (I-C),
\quad n \geq 1.
\eqno (22)
$$
In other words, all what we can say is that $\b F_n$ is what we
could call a left quasi-orthogonal matrix polynomial of order $n$
with respect to $d\Omega$, that is, a non null matrix polynomial
with $\deg \b F_n \leq n$ and left orthogonal to $I, Ix, Ix^2,
\dots, Ix^{n-2}$ with respect to the measure $d\Omega$ (see for
example [5] or [9] for introducing the analogous conception in the
scalar case).

Notice that, when the measure $d\mu$ is symmetric, both $\b F_n$ and $d\Omega$ are
diagonal. Then, our quasi-orthogonal matrix polynomials provide two sequences of
scalar OP on $[-1,1]$, that is, we recover Szeg\H o's result.

The complex recurrence formula for VSOF provides two real recurrence relations for the
corresponding matrix polynomials.

\noindent {\bf Proposition 5.}
{\sl The matrix polynomials associated to $d\mu$ satisfy the recurrence relations
$$
\eqalign{
& x \b F_n(x) = \b F_{n+1}(x) + L_n \b F_n(x) + M_n \b F_{n-1}(x), \quad n \geq 1, \cr
& \b F_n(x) Y(x) = J \b F_{n+1}(x) +
\widetilde L_n \b F_n(x) + \widetilde M_n \b F_{n-1}(x), \quad n \geq 1, \cr
}
\eqno (23)
$$
where
$$
\eqalign{
& Y(x) = \pmatrix{0 & 1 \cr 1-x^2 & 0}, \cr
& M_n =
{1 \over 4} \det (I - H_{2n-1}) (I - H_{2n}) (I + H_{2n-2}), \cr
& \widetilde M_n =
- {1 \over 4} \det (I - H_{2n-1}) (I - H_{2n}) J (I + H_{2n-2}), \cr
& L_n =
{1 \over 2} \left\{ (I - H_{2n}) H_{2n-1} - H_{2n+1} (I + H_{2n}) \right\},
\cr & \widetilde L_n =
{1 \over 2} \left\{ (I - H_{2n}) H_{2n-1} J - J H_{2n+1} (I + H_{2n}) \right\}, \cr
}
\eqno (24)
$$
and $H_n$ are the Schur matrices related to $d\mu$.
}

\noindent {\bf Proof.}
VSOF, like SOF, satisfy $\overline{\b f}_n(z^{-1})=\b f_n(z)$.
Therefore, from the recurrence relation in Proposition 3 we get
$$
z^{-1} \b f_n(z) =  (I-iJ) \b f_{n+1}(z) +
\overline{\b L}_n \b f_n(z) + \overline{\b M}_n \b f_{n-1}(z).
$$
Appropriate linear combinations of this new recurrence relation and the original one
give
$$
\eqalign{
& x \b f_n(z) = \b f_{n+1}(z) + \Re \b L_n \b f_n(z) + \Re \b M_n \b f_{n-1}(z), \cr
& y \b f_n(z) = J \b f_{n+1}(z) + \Im \b L_n \b f_n(z) + \Im \b M_n \b f_{n-1}(z). \cr
}
\eqno (25)
$$
From (17) we get
$$
y \b f_n(z) = \b F_n(x) Y(x) \pmatrix{1 \cr y \cr},
\quad Y(x) = \pmatrix{0 & 1 \cr 1-x^2 & 0}.
\eqno (26)
$$
Introducing (17) and (26) in (25) we find two relations for $\b F_n$ with the
form
$$
A(x) \pmatrix{1 \cr y} = B(x) \pmatrix{1 \cr y}, \quad A(x),B(x) \in {\cal P}^{(2,2)},
$$
that are true if $x=(z+z^{-1})/2, y=(z-z^{-1})/2i$ for any $z \neq 0$.
Evaluating in $z$ and $z^{-1}$ we obtain
$$
A(x) \pmatrix{1 & 1 \cr y & -y} = B(x) \pmatrix{1 & 1 \cr y & -y}
$$
and, therefore, it must be $A(x)=B(x)$. Taking into
account the expressions for $\b M_n$ and $\b L_n$ given in Proposition 3, we see that
the two equalities that we find in this way are exactly the desired recurrence
relations for $\b F_n$.
$\qed$

From Proposition 4 we see that
$$
\eqalign{
& \b F_0(x) = C, \cr
& \b F_1(x) = (xI - I + H_1)C + I. \cr
}
\eqno (27)
$$
Therefore, starting from the Schur matrices $H_n$ associated to
$d\mu$, the first recurrence relation in (23), together with the
expressions (27) for the two first matrix polynomials, let us
obtain the complete sequence of matrix polynomials associated to
$d\mu$. Conversely, suppose that we have an arbitrary sequence
$(H_n)_{n\geq 0}$ of $2 \times 2$ real symmetric traceless
matrices with $H_0 = 2 C - I$ and $|\det H_n| < 1$ for $n \geq 1$.
If a sequence $({\b F_n})_{n\geq  0}$ of matrix polynomials
satisfies a recurrence relation like the first one in (23) with
$M_n$ and $L_n$ given by (24) and the initial conditions (27),
then the matrix polynomials $({\b F_n})_{n\geq  0}$ are associated
to some measure on $\T$ (and, therefore, they are left
quasi-orthogonal with respect some matrix measure on $[-1,1]$). To
see this, just notice that the matrix sequence $(H_n)_{n\geq 0}$
provides, through the relation $a_n = H_n^{(11)} + i H_n^{(12)}$,
a complex sequence $(a_n)_{n\geq 0}$ such that $a_0 = 1$ and
$|a_n| < 1$ for $n \geq 1$ ($A^{(kj)}$ denotes the $(k,j)$-th
element of the matrix $A$). Now, it is well known that this
conditions for $a_n$ ensure that the complex polynomials
$(\phi_n)_{n\geq0}$ defined by (1) form a sequence of monic OP
with respect to some positive measure $d\mu$ on $\T$. We have
shown that the measure $d\mu$ generates an associated matrix
measure $d\Omega$ on $[-1,1]$ and that the OP $(\phi_n)_{n\geq0}$
let us construct a sequence of left quasi-orthogonal matrix
polynomials with respect to $d\Omega$ satisfying (23), (24), and
(27). This sequence must be $({\b F_n})_{n\geq 0}$.

\vskip 0.5cm

 \noindent {\bf 4. Semi-orthogonal functions and left orthogonal matrix
polynomials}

\vskip 0.25cm

We have discovered that the generalization of Szeg\H o's method leads in
general, not to a sequence of scalar OP, neither a sequence of LOMP, but to a
sequence $({\b F_n})_{n\geq 0}$ of left quasi-orthogonal matrix polynomials
with respect to the matrix measure $d\Omega$. However, if $d\Omega$ is a
positive matrix measure then there exists a sequence of LOMP with respect to
$d\Omega$ iff [6],
$$
\int_{-1}^1 p^T(x) d\Omega(x) p(x) \neq 0,
\quad \forall p \in \C^{(2,1)}[x] \backslash \{0\}
\eqno (28)
$$
As we see in the following proposition, when a matrix measure is associated to
a measure on $\T$, it is possible to give simple conditions equivalent to
(28).

\noindent {\bf Proposition 6.}
{\sl Let $d\mu$ be an arbitrary measure on $\T$ and let $d\Omega$ be the matrix
measure associated to $d\mu$. Then, $d\Omega$ is positive iff $d\mu$ is
positive. Moreover, when $d\mu$ is positive the following statements are
equivalent:
\item{i)} There exists a sequence of LOMP with respect to $d\Omega$.
\item{ii)} There exists a sequence of OP with respect to $d\mu$.
\item{iii)} supp$(d\mu)$ is infinite.
\item{iv)} supp$(d\Omega)$ is infinite.
}

\noindent {\bf Proof.} Let us suppose that $d\mu$ is positive.
Then $\vert d\sigma \vert \leq d\rho$ and, thus, $d\rho$ is
positive (see Remark 6). In order to prove the positivity of
$d\Omega$ we have just to see that $\int_a^b d\Omega(x)$ is a
nonnegative definite matrix for all $a,b \in [-1,1]$, which is
equivalent to say that its trace and determinant are both
nonnegative. Since $d\rho$ is positive,
$$
{\rm tr} \int_a^b d\Omega(x) = \int_a^b (2-x^2) d\rho(x) \geq 0.
$$
Moreover, taking into account that $\vert d\sigma \vert \leq d\rho$ we get that
$$
\eqalign{
{\rm det} \int_a^b d\Omega(x) & =
\int_a^b d\rho(x) \int_a^b (1-x^2) d\rho(x) -
\left( \int_a^b \sqrt{1-x^2} d\sigma(x) \right)^2 \geq \cr
& \geq \int_a^b d\rho(x) \int_a^b (1-x^2) d\rho(x) -
\left( \int_a^b \sqrt{1-x^2} d\rho(x) \right)^2 \geq 0,
}
$$
where we have used the Schwarz's inequality. Thus, if $d\mu$ is
positive then $d\Omega$ is positive too.

To see the converse first notice that if $p \in \C^{(2,1)}[x]$
then $f(z)=\pmatrix{1&y}p(x)$ is a Laurent polynomial. Even more,
for every $f \in \Lambda$ there is a unique decomposition
$f(z)=\pmatrix{1&y}p(x), p \in \C^{(2,1)}[x]$. This decomposition
holds iff $p^T=(p_1,p_2)$ with $p_1(x)=(f(z)+f(z^{-1}))/2$ and
$p_2(x)=(f(z)-f(z^{-1}))/2y$. Now, using the Tchebychev
polynomials of first and second kind we see that $p_1,p_2 \in
\C[x]$. This provides an isomorphism between $\Lambda$ and
$\C^{(2,1)}[x]$. Let us consider an arbitrary $f \in \Lambda$ and
the corresponding $p \in \C^{(2,1)}[x]$. Then,
$$
\int_0^{2\pi} \vert f(e^{i\theta}) \vert^2 d\mu(\theta) =
2\int_{-1}^1 p^T(x) d\Omega(x) p(x).
$$
If $d\Omega$ is a positive matrix measure then
$\int_{-1}^1 p^T(x) d\Omega(x) p(x) \geq 0$
for all $p \in \C^{(2,1)}[x]$.
Therefore,
$\int_0^{2\pi} \vert f(e^{i\theta}) \vert^2 d\mu(\theta) \geq 0$
for all $f \in \Lambda$
and, thus, $d\mu$ is po\-si\-ti\-ve too.

Now, assume that $d\mu$ is positive. The equivalence between ii)
and iii) is known. From (28) and above results we see that the
statement i) means that $\int_0^{2\pi} \vert f(e^{i\theta})
\vert^2 d\mu(\theta) \neq 0$ for all $f \in \Lambda \backslash
\{0\}$, which holds iff $d\mu$ has an infinite support. So it is
proved that i) is equivalent to iii). The equivalence between iii)
and iv) is just a consequence of the following facts that are true
for any positive measure $d\mu$ (see Remark 6): supp$(d\rho)$ is
infinite iff supp$(d\mu)$ so is; supp$(d\sigma) \subset$
supp$(d\rho)$ and, hence, supp$(d\Omega)=$ supp$(d\rho)$.

In what follows we will suppose again that $d\mu$ is a positive measure on $\T$
with infinitely many points in the support. Then, the next proposition gives a
sequence of LOMP with respect to the related matrix measure in terms of the
associated matrix polynomials.

\noindent {\bf Proposition 7.}
{\sl Let $d\Omega$ and $(\b F_n)_{n\geq0}$ be the matrix measure and the matrix
polynomials  associated to $d\mu,$ respectively. Then, the matrix polynomials
$(\b P_n)_{n\geq0}$ given by
$$
\eqalign{
& \b P_n(x) = \alpha_n \b F_{n+1}(x) + \beta_n \b F_n(x), \quad n \geq 0, \cr
& \alpha_n = I-C, \quad \beta_n = \pmatrix{1 & r_n \cr 0 & 0 \cr}, \cr
& r_n = \Im a_{2n} / (1+\Re a_{2n}), \cr
}
$$
define a sequence of LOMP with respect $d\Omega$, where $a_n$ are the Schur parameters
related to $d\mu$. Moreover,
$$
\eqalign{
& \b F_n(x) = \widetilde \alpha_n \b P_n(x) + \widetilde \beta_n \b P_{n-1}(x),
\quad n \geq 0, \cr
& \widetilde \alpha_n = C,
\quad \widetilde \beta_n = \pmatrix{0 & -r_n \cr 0 & 1 \cr}, \cr
}
$$
with the convention $\b P_{-1}=0$.
}

\noindent{\bf Proof.}
Let $({\b P_n})_{n\geq 0}$ be an arbitrary sequence of LOMP with respect to $d\Omega$. From
the algebraic and orthogonality properties of $\b F_n$ and $\b P_n$ we have that, for
$n \geq 0$
$$
\eqalign{
& \b P_n(x) = \alpha_n \b F_{n+1}(x) + \beta_n \b F_n(x),
\quad \alpha_n, \beta_n \in \R^{(2,2)}, \cr
& \b F_n(x) = \widetilde \alpha_n \b P_n(x) + \widetilde \beta_n \b P_{n-1}(x),
\quad \widetilde \alpha_n, \widetilde \beta_n \in \R^{(2,2)}, \cr
}
\eqno (29)
$$
where $\b P_{-1}=0$ and the matrix coefficients must satisfy the relations
$$
\eqalign{
& \alpha_n \widetilde \alpha_{n+1} = \beta_n \widetilde \beta_n = 0, \cr
& \alpha_n \widetilde \beta_{n+1} + \beta_n \widetilde \alpha_n = I.
}
\eqno (30)
$$
Now, we proceed to determine $\alpha_n$, $\beta_n$, $\widetilde
\alpha_n$, $\widetilde \beta_n$ by imposing on $\b P_n$ the
conditions (I) and (II) given after Remark 7.

From (29) and Proposition 4 i) we get
$$
\b P_n(x) = \alpha_n C x^{n+1} +
\left\{ \alpha_n \pmatrix{\eta_n & 0 \cr \gamma_n & 1} + \beta_n C \right\} x^n +
\dots, \quad n \geq 0,
$$
where the dots mean terms with degree less than $n$. Therefore, (I) is equivalent to
$$
\eqalign{
& \alpha_n C = 0, \quad n \geq 0, \cr
& \alpha_n \pmatrix{\eta_n & 0 \cr \gamma_n & 1} +
\beta_n C \quad {\rm nonsingular}, \quad n \geq 0. \cr
}
\eqno (31)
$$

The quasi-orthogonality of $\b F_n$ implies that, for all
$\alpha_n, \beta_n \in \R^{(2,2)}$, the matrix polynomial $\b P_n$ given in (29)
is orthogonal to $I x^k$, $0 \leq k \leq n-2$, with respect to $d\Omega$. For
$I x^{n-1}$, we can use (21) and (22) to obtain
$$
\int_{-1}^1 \b P_n(x) \, d\Omega(x) \, x^{n-1} = {1 \over 2} \beta_n C_n (I-C),
\quad n \geq 1.
$$
Since this integral must vanish, with the aid of (20) we find
$$
\int_{-1}^1 \b P_n(x) \, d\Omega(x) \, x^n =
{1 \over 2} \left\{ \alpha_n  C_{n+1} (I-C) + \beta_n C_n \right\} (\Gamma_n^{-1})^T,
\quad n \geq 0.
$$
Hence, (II) is equivalent to
$$
\eqalign{
& \beta_n C_n (I-C) = 0, \quad n \geq 1, \cr
& \alpha_n  C_{n+1} (I-C) + \beta_n C_n
\quad {\rm nonsingular}, \quad n \geq 0. \cr
}
\eqno (32)
$$

If $V_n \in \R^{(2,2)}$ is nonsingular for all $n \geq 0$, then
$$
\eqalign{
& \alpha_n = V_n (I-C), \quad n \geq 0, \cr
& \beta_0 = V_0 C, \cr
& \beta_n = V_n C C_n^{-1}, \quad n \geq 1, \cr
}
$$
are solutions of (31) and (32). The expressions given in
the proposition for $\alpha_n$ and $\beta_n$ correspond to the choice $V_0=I$ and
$$
V_n = \pmatrix{{1 \over (C_n^{-1})^{(11)}} & 0 \cr 0 & 1}, \quad n \geq 1.
$$
The relations (30) give then $\widetilde \alpha_n$ and $\widetilde \beta_n$.
$\qed$

\noindent {\bf Remark 8.}
Taking into account Proposition 4 and Proposition 7, we get that the leading coefficient
of $\b P_n$ is $\Gamma_n$ (see (20)) and
$$
\int_{-1}^1 \b P_n(x) \, d\Omega(x) \, \b P_n(x) = 2^{-2n}
\pmatrix{\varepsilon_{2n} (1+\Re a_{2n})^{-1} & 0 \cr
0 & {1 \over 4} \varepsilon_{2n+1} (1+\Re a_{2n+2}) \cr},
$$
where $\varepsilon_n$ is given in (5).
Therefore, for the matrix measure $d\Omega$, the monic LOMP are
$\widetilde{\b P}_n = \Gamma_n^{-1} \b P_n$ while
$\hat \b P_n = W_n \b P_n$ are left orthonormal polynomials (LONP), being
$$
W_n = 2^n
\pmatrix{\kappa_{2n} (1+\Re a_{2n})^{1/2} & 0 \cr
0 & 2 \kappa_{2n+1} (1+\Re a_{2n+2})^{-1/2} \cr}.
\eqno (33)
$$
with $\kappa_n=\varepsilon_n^{-1/2}$.

 Given a matrix measure, LOMP are determined up to multiplication on the left by a
nonsingular constant matrix. Therefore, monic LOMP are unique but LONP are
defined up to multiplication on the left by an orthogonal constant matrix. Thus, LONP
can be fixed if we ask for their coefficients to be symmetric and positive definite.
We will refer to the standard LONP when this choice is made.
As for the measure $d\Omega$, we see that $(\hat \b P_n)_{n\geq0}$ is not the sequence
of standard LONP because the corresponding leading coefficients $W_n \Gamma_n$ are not
symmetric. The following proposition gives the standard LONP in this case.

\noindent {\bf Proposition 8.}
{\sl Let $d\Omega$ be the matrix measure associated to $d\mu$ and let
$(\b P_n)_{n\geq0}$ be the corresponding LOMP defined in Proposition 7. Then, the
sequence of standard LONP $(\b Q_n)_{n\geq0}$ with respect to $d\Omega$ is given by
$\b Q_n = \Xi_n^T W_n \b P_n$, where
$$
\Xi_n = {1 \over \sqrt{\det \Theta_n}} \Theta_n,
\quad \Theta_n = K_n - J K_n J = K_n + {\rm adj} K_n^T,
$$
being $K_n = W_n \Gamma_n$ with $\Gamma_n$ and $W_n$ defined in (20) and (33)
respectively.
}

\noindent {\bf Proof.}
Since $\hat \b P_n = W_n \b P_n$ is a LONP, $\b Q_n$ will be a LONP too iff $\Xi_n$ is
an orthogonal matrix. To see that $\Xi_n$ is indeed orthogonal, first notice that
$JAJ=-{\rm adj}A^T$ for all $A \in \R^{(2,2)}$. Therefore, the nonnegative definite
matrix
$(A - J A J)^T (A - J A J) =
A^T A + {\rm adj} (A^T A) + A^T ({\rm adj} A^T) + ({\rm adj} A) A =
({\rm tr}(A^T A) + 2 (\det A)) I$
is a multiple of the identity. Taking determinants in above expression we
see that this multiple is the nonnegative factor $\det (A + {\rm adj} A^T)$.
For $A=K_n,$ this factor can not vanish because $\det K_n = \det W_n > 0$. Thus,
$\Xi_n^T \Xi_n=I$.

 Now, the leading coefficient of $\hat \b P_n$ is $K_n = W_n \Gamma_n$. So, the leading
coefficient of $\b Q_n$ is $\Xi_n^T K_n = K_n^T K_n + (\det K_n) I$, which is symmetric
and positive definite because $\det K_n > 0$. Hence, $(\b Q_n)_{n\geq0}$ are the
standard LONP.
$\qed$

In the next section we will deal with the strong asymptotics of
LOMP with respect to a matrix measure associated to a measure on
$\T$. As usual, we will take the standard LONP as a reference to
express the asymptotic behavior.

\vskip 0.5cm

\noindent {\bf 5. Asymptotics of semi-orthogonal functions and matrix OP}

\vskip 0.25cm

Once we have shown above connection between scalar and matrix OP,
it is natural to take advantage of known properties for OP on $\T$
to develop new results about the more unfamiliar world of matrix
OP. As an example we present here the implications of the
asymptotics of OP on $\T$ when Szeg\H o's condition,
$$
\int_0^{2\pi} \log \mu'(\theta) \; d\theta > -\infty,
\eqno (34)
$$
for the measure $d\mu$ holds (as it is usual, $\mu'$ denotes the
Radon-Nikodym derivative of the absolutely continuous part of
$d\mu$ with respect to the Lebesgue measure $d\theta$). It can be
proved [9, 24] that Szeg\H o's condition is equivalent to
$(a_n)_{n\geq 0} \in \ell^2$, which, in sight  of (5), means that
$\varepsilon = \lim_n \varepsilon_n > 0$ (or, in other words,
$\kappa = \lim_n \kappa_n < \infty$). Thus, a necessary condition
for (34) is $\lim_n a_n = 0$.

When Szeg\H o's condition holds, asymptotic properties of OP on $\T$ are given in
terms of the function
$$
D(d\mu;z) = \exp \left( {1 \over 4\pi} \int_0^{2\pi} \log \mu'(\theta)
\, {e^{i\theta} + z \over e^{i\theta} - z} \; d\theta \right), \quad |z| \neq 1,
$$
which, for $|z| < 1$, it is known as Szeg\H o's function for the measure $d\mu$, and
satisfies the following remarkable property [9, 24],
$$
\lim_{r \to 1^{-}} D(d\mu;re^{i\theta}) \overline{D(d\mu;re^{i\theta})} = \mu'(\theta),
\quad {\rm a.e.}.
$$
Notice that, for $d\widetilde\mu$, the symmetric measure of $d\mu$,
$$
D(d\widetilde\mu;z) = \overline{D(d\mu;\overline{z})} = D(d\mu;z^{-1})^{-1}.
$$
It can be proved [9, 24] that, under (34),
$$
\kappa = {1 \over \sqrt{2\pi}} D(d\mu;0)^{-1}
\eqno (35)
$$
and the orthonormal polynomials satisfy
$$
\eqalign{
& \lim_n \varphi_n(z) = 0, \quad |z| < 1, \cr
& \lim_n \varphi_n^*(z) = {1 \over \sqrt{2\pi}} D(d\mu;z)^{-1}, \quad |z| < 1, \cr
}
$$
where the convergence is uniform on compact sets. Therefore,
$$
\varepsilon = 2\pi D(d\mu;0)^2
$$
and for the monic OP we get
$$
\eqalign{
& \lim_n \phi_n(z) = 0, \quad |z| < 1, \cr
& \lim_n \phi_n^*(z) = D(d\mu;0) D(d\mu;z)^{-1}, \quad |z| < 1. \cr
}
\eqno (36)
$$
From above results the following asymptotics of VSOF follows straightforward
$$
\eqalign{
\lim_n 2^n z^n \b f_n(z)
&= \pmatrix{1 \cr i} D(d\mu;0) D(d\mu;z)^{-1}, \quad 0 < |z| < 1, \cr
\lim_n 2^n z^{-n} \b f_n(z)
&= \pmatrix{1 \cr -i} D(d\mu;0) D(d\mu;z), \quad |z| > 1, \cr
}
\eqno (37)
$$
where, again, the convergence is uniform on compact sets.

  As for matrix OP on $[-1,1]$, some general results are known, but they are not so good
as previous ones. More precisely, let us suppose that a positive matrix measure
$d\omega$ on $[-1,1]$ satisfies the Szeg\H o's matrix condition
$$
\int_{-1}^1 \log \det \omega'(x) \, {dx \over \sqrt{1-x^2}} > -\infty,
$$
where $\omega'$ is the Radon-Nikodym derivative of the absolutely continuous part of
$d\omega$ with respect to the Lebesgue scalar measure $d\theta$.
If $(\b q_n)_{n\geq0}$ is the sequence of standard LONP with respect to
$d\omega$ then, for any other sequence of LONP $(\b p_n)_{n\geq0}$ such that
$\b p_n = \xi_n \b q_n$ with $\sn \xi 0$ convergent, we have that
$$
\lim_n {z^n \b p_n(x)} = {1 \over \sqrt{2\pi}} \b D(d\omega;z)^{-1},
\quad |z| < 1,
\eqno (38)
$$
where the convergence is uniform on compact sets and $\b
D(d\omega;z)$ is certain matrix-valued analytic function on $|z| <
1$ without zeros there [2]. The Szeg\H o's matrix function $\b
D(d\omega;z)$ is uniquely determined by $\omega'$ and satisfies
the boundary condition
$$
\lim_{r \to 1^{-}} \b D(d\omega;re^{i\theta}) \b D(d\omega;re^{i\theta})^* =
\omega'(\cos\theta)|\sin\theta|, \quad {\rm a.e.}.
$$

  Unfortunately, an explicit expression for $\b D(d\omega;z)$ in terms of $\omega'$ is
not a\-vai\-la\-ble. On that score, all what we can state is that [2]
$$
\b D(d\omega;z) = \buildrel\curvearrowright\over{\int_0^{2\pi}}
\exp \left( \b M(\theta)  {e^{i\theta} + z \over e^{i\theta} - z} \; d\theta \right)
\zeta, \quad |z| < 1,
$$
where $\zeta$ is a constant orthogonal matrix factor depending on $\lim_n \xi_n$, and
$\b M(\theta)$ is a Hermitian matrix-valued integrable function on $[0,2\pi)$ such that
$$
{\rm tr} \b M(\theta) = \log \det \{ \omega'(\cos\theta)|\sin\theta| \},
\quad \theta \in [0,2\pi).
$$
The symbol $\buildrel\curvearrowright\over{\int_0^{2\pi}}$ means the multiplicative
integral
$$
\buildrel\curvearrowright\over{\int_0^{2\pi}} \exp(F(\theta)) \, d\theta =
\lim_n \prod_{k=1}^n \exp(F(t_k)) (\theta_k - \theta_{k-1}),
$$
where $t_k \in [\theta_{k-1},\theta_k)$ and the limit is taken in the usual sense
over the partitions
$0=\theta_0 < \theta_1 < \theta_2 < \cdots < \theta_{n-1} < \theta_n=2\pi$ of the
interval $[0,2\pi)$.

When the matrix measure is associated to a measure on $\T$ above results can
be improved by translating the better known asymptotics of OP on $\T$ to matrix OP on
$[-1,1]$. As a consequence we can obtain in this case an explicit expression for the
Szeg\H o's matrix function.

\noindent {\bf Theorem 1.}
{\sl Let
$$
d\Omega(x) = {1 \over 2}
\pmatrix{d\rho(x) & \sqrt{1-x^2} d\sigma(x) \cr
\sqrt{1-x^2} d\sigma(x) & (1-x^2) d\rho(x) \cr}
$$
be a positive matrix measure on $[-1,1]$ ($d\rho$ and $d\sigma$ are scalar measures on
$[-1,1]$) that satisfies the condition
$$
\int_{-1}^1 \log \det \Omega'(x) \, {dx \over \sqrt{1-x^2}} > -\infty.
$$
Let ${\cal R}(z),{\cal I}(z), \gamma$ be
$$
\eqalign{
{\cal R}(z)
%& = {1 \over 4\pi} \int_0^{\pi} \log \det \Omega'(\cos\theta)
%\; \Re_\theta \!\! \left( {e^{i\theta}+z \over e^{i\theta}-z} \right) d\theta \cr
& = {1-z^2 \over 4\pi} \int_{-1}^1
{\log \det \Omega'(x) \over \sqrt{1-x^2}} \;
{dx \over 1-2xz+z^2}, \quad |z| \neq 1, \cr
{\cal I}(z)
%& = {1 \over 4\pi} \int_0^{\pi}
%\log \left\{ {\rho'(\cos\theta) + \sigma'(\cos\theta) \over
%\rho'(\cos\theta) - \sigma'(\cos\theta)} \right\}
%\Im_\theta \!\! \left( {e^{i\theta}+z \over e^{i\theta}-z} \right) d\theta  \cr
& = -{z \over 2\pi} \int_{-1}^1
\log \left\{ {\rho'(x) + \sigma'(x) \over \rho'(x) - \sigma'(x)} \right\}
{dx \over 1-2xz+z^2}, \quad |z| \neq 1, \cr
\gamma
& = -{1 \over 4\pi} \int_{-1}^1
\log \left\{ {\rho'(x) + \sigma'(x) \over \rho'(x) - \sigma'(x)} \right\} \, dx.
}
$$
%where $\Re_\theta$ and $\Im_\theta$ are respectively the real and imaginary part
%operators with conjugation acting only on $e^{i\theta}$-dependence.
and let $(\b Q_n)_{n\geq0}$ be the sequence of standard LONP with respect to $d\Omega$.
Then, for $x \in \C \setminus [-1,1]$, if we write $z = x + \sqrt{x^2-1}$ with the
choice of $\sqrt{x^2-1}$ such that $|z|<1$, we have that
$$
\lim_n z^n \b Q_n(x) = {1 \over \sqrt{2\pi}} \b D(d\Omega;z)^{-1}
$$
where
$$
\b D(d\Omega;z) =
{1 \over \sqrt{9 + 4 \gamma^2}} \pmatrix{1 & 0 \cr 0 & -\sqrt{x^2-1}}
\exp(I {\cal R}(z) + J {\cal I}(z)) \pmatrix{3 & -2\gamma \cr 2\gamma z & 3z},
$$
being the convergence uniform on compact sets.
}

\noindent {\bf Proof.}
Notice first that Szeg\H o's matrix condition implies that $d\Omega$ has an
infinite support. Thus, from Proposition 6 we see that there exist LOMP with respect
to $d\Omega$ and that the matrix measure $d\Omega$ is associated to some positive
measure $d\mu$ on $\T$ with an infinite support (therefore, there exist OP with respect
to $d\mu$).

The expression that gives the Szeg\H o condition for $d\mu$ can be rewritten in the
following way
$$
\eqalign{
\int_0^{2\pi} \log \mu'(\theta) \; d\theta &=
\int_{-1}^1 \left[\log \left( \nu'_1(x) \sqrt{1-x^2} \right) +
\log \left( \nu'_2(x) \sqrt{1-x^2} \right) \right] {dx \over \sqrt{1-x^2}} \cr
&= \int_{-1}^1 \log \left\{ {1 \over 4} \left(\rho'(x)^2 - \sigma'(x)^2 \right) (1-x^2)
\right\} {dx \over \sqrt{1-x^2}} \cr
&= \int_{-1}^1 \log \det \Omega'(x) \, {dx \over \sqrt{1-x^2}}. \cr
}
$$
Thus, Szeg\H o's matrix condition for $d\Omega$ is equivalent to
Szeg\H o's condition for $d\mu$. Hence, under the assumptions of
the theorem, the Szeg\H o's function $D(d\mu;z)$ go\-verns the
asymptotic behavior of the VSOF related to $d\mu$ in the way shown
in (37).

 Then, from (18) we find for the quasi-orthogonal matrix polynomials associated to
$d\mu$ that
$$
\lim_n 2^n z^{-n} \b F_n(x) = D(d\mu;0) \; \D(d\mu;z) \pmatrix {1 & 0 \cr 0 & iy}^{-1},
\quad |z| > 1,
$$
where
$$
\D(d\mu;z) = \pmatrix{D_s(d\mu;z) & iD_a(d\mu;z) \cr -iD_a(d\mu;z) & D_s(d\mu;z)}
= I D_s(z) + iJ D_a(z),
$$
and $D_s(d\mu;z)$, $D_a(d\mu;z)$ are what we could call the symmetric and
antisymmetric part of $D(d\mu;z)$, that is,
$$
\eqalign{
& D_s(z) = {1 \over 2} (D(d\mu;z) + D(d\widetilde\mu;z)), \cr
& D_a(z) = {1 \over 2} (D(d\mu;z) - D(d\widetilde\mu;z)). \cr
}
$$
Notice that, when $d\mu$ is symmetric, $D_s(d\mu;z)=D(d\mu;z)$ and
$D_a(d\mu;z)=0$. Hence, the matrix $\D(d\mu;z)$ is diagonal. This is natural,
because in this case the quasi-orthogonal polynomials are diagonal.

Let us define
$$
\eqalign{
& {\cal R}(z) = {1 \over 4\pi} \int_0^{2\pi} \log \mu'(\theta)
\; \Re_\theta \!\! \left( {e^{i\theta}+z \over e^{i\theta}-z} \right) d\theta, \cr
& {\cal I}(z) = {1 \over 4\pi} \int_0^{2\pi} \log \mu'(\theta)
\; \Im_\theta \!\! \left( {e^{i\theta}+z \over e^{i\theta}-z} \right) d\theta, \cr
}
$$
where $\Re_\theta$ and $\Im_\theta$ are real and imaginary part
operators with conjugation acting only on
$e^{i\theta}$-dependence. Then, $D(d\mu;z) = \exp({\cal
R}(z)+i{\cal I}(z))$ and
$$
\D(d\mu;z) = \exp{\cal R}(z)
\pmatrix{\cos{\cal I}(z) & -\sin{\cal I}(z) \cr \sin{\cal I}(z) & \cos{\cal I}(z)}
= \exp(I {\cal R}(z) - J {\cal I}(z)).
$$
Notice that
$$
\D(d\widetilde\mu;z) = \exp(I {\cal R}(z) + J {\cal I}(z)) = \D(d\mu;z)^T
= \overline{\D(d\mu;\overline{z})} = \D(d\mu;z^{-1})^{-1}.
$$

The asymptotic behavior of the LOMP given in Proposition 7 can now
be deduced. Since $\lim_n a_n = 0$ we get
$$
\lim_n 2^n z^{-n} \b P_n(x)
= D(d\mu;0) \pmatrix{1 & 0 \cr 0 & z/2} \D(d\mu;z)
  \pmatrix {1 & 0 \cr 0 & iy}^{-1}, \quad |z| > 1.
$$

As for the standard LONP, Proposition 6 together with (20), (33) and (35) give
$$
\lim_n z^{-n} \b Q_n(x) =  {1 \over \sqrt{2\pi}} {1 \over \sqrt{9+4\gamma^2}}
\pmatrix{3 & 2\gamma z \cr -2\gamma & 3z} \D(d\mu;z)
\pmatrix {1 & 0 \cr 0 & iy}^{-1}, \quad |z| > 1,
$$
where $\gamma = \lim_n \gamma_n$. This limit exits because, from the expression for
$\gamma_n$ given in Proposition 4 i) and relation (3), we see that
$$
\gamma = {1 \over 2} \Im \left(\lim_n b_n\right)
= {1 \over 2} \Im \left( \sum_{k=1}^\infty a_k \overline{a}_{k-1} \right),
$$
and the convergence of $\sum_{k=0}^\infty a_{k+1} \overline{a}_k$ follows from the
convergence of $\sum_{k=0}^\infty |a_k|^2$ and the Schwarz's inequality.

 By comparing with (38) we see that, in our case, the asymptotics of
$(\b Q_n)_{n\geq0}$ is governed by the Szeg\H o's matrix function
$$
\b D(d\Omega;z) = {1 \over \sqrt{9+4\gamma^2}} \pmatrix {1 & 0 \cr 0 & -iy}
\D(d\mu;z)^T
\pmatrix{3 & -2\gamma \cr 2\gamma z & 3z}, \quad |z| < 1.
$$

To complete the proof it only remains to see that ${\cal R}(z)$, ${\cal I}(z)$
and $\gamma$ are given by the expressions that appear in the theorem. We can rewrite
${\cal R}(z)$ and ${\cal I}(z)$ in terms of the matrix measure $d\Omega$ in the
following way
$$
\eqalign{
{\cal R}(z) &= {1 \over 4\pi} \int_0^{\pi}
\log \{ \nu'_1(\cos\theta) \nu'_2(\cos\theta) \sin^2\theta \}
\; \Re_\theta \!\! \left( {e^{i\theta}+z \over e^{i\theta}-z} \right) d\theta \cr
&= {1 \over 4\pi} \int_0^{\pi}
\log \det \Omega'(\cos\theta)
\; \Re_\theta \!\! \left( {e^{i\theta}+z \over e^{i\theta}-z} \right) d\theta, \cr
{\cal I}(z) &= {1 \over 4\pi} \int_0^{\pi}
\log \{ \nu'_1(\cos\theta) / \nu'_2(\cos\theta) \}
\; \Im_\theta \!\! \left( {e^{i\theta}+z \over e^{i\theta}-z} \right) d\theta \cr
&= {1 \over 4\pi} \int_0^{\pi}
\log \left\{ {\rho'(\cos\theta) + \sigma'(\cos\theta) \over
\rho'(\cos\theta) - \sigma'(\cos\theta)} \right\}
\Im_\theta \!\! \left( {e^{i\theta}+z \over e^{i\theta}-z} \right) d\theta. \cr
}
$$
Now, the change of variables $x=\cos\theta$ gives the desired expressions for
${\cal R}(z)$ and ${\cal I}(z)$.

Besides, we can give an expression for $\gamma = {1 \over 2} \Im (\lim_n b_n)$ in terms
of $d\Omega$. Notice first that $b_n = \overline{{\phi_n^*}'(0)}$. Since
$\sn {\phi^*} 0$ is a sequence of analytic functions in the complex plane that
converges uniformly on compact subsets of $|z|<1$, we can write
$\lim_n {\phi_n^*}'(z) = (\lim_n {\phi_n^*})'(z)$ for $|z|<1$ (see [22]). Thus, from
(36) we get
$$
\lim_n b_n = -\overline{(\log D)'(d\mu;0)}
= -{1 \over 2\pi} \int_0^{2\pi} \log \mu'(\theta) \; e^{i\theta} \, d\theta.
$$
Therefore,
$$
\eqalign{
\gamma
&= -{1 \over 4\pi} \int_0^{2\pi} \log \mu'(\theta) \, \sin\theta \, d\theta \cr
&= -{1 \over 4\pi} \int_0^{\pi} \log \{ \nu'_1(\cos\theta)/\nu'_2(\cos\theta) \}
\, \sin\theta \, d\theta \cr
&= -{1 \over 4\pi} \int_{-1}^1
\log \left\{ {\rho'(x) + \sigma'(x) \over \rho'(x) - \sigma'(x)} \right\} \, dx, \cr
}
$$
which completes the proof.
$\qed$

\vskip 0.5cm

\noindent {\bf Acknowledgements}

\vskip 0.25cm

This work was supported by Direcci\'on General de ense\~nanza Superior (DGES) of Spain
under grant PB 98-1615.

\vskip 0.5cm

\noindent {\bf References}

\vskip 0.25cm

\item{[1]}
M. Alfaro, M.J. Cantero, L. Moral,
Semi-orthogonal functions and or\-tho\-go\-nal polynomials on the unit circle,
{\sl J. Comput. Appl. Math.} {\bf 99} (1998) 3--14.

\item{[2]}
A.I. Aptekarev, E.M. Nikishin,
The scattering problem for a discrete Sturm-Liouville operator,
{\sl Math. USSR Sb.} {\bf 49} (1984) 325--355.

\item{[3]}
Yu.W. Berezanskii,
{\sl Expansions in Eigenfunctions of Self-adjoint Operators},
Transl. Math. Monographs 17, Amer.Math.Soc., 1968.

\item{[4]}
M.J. Cantero,
{\sl Polinomios ortogonales sobre la circunferencia unidad. Modificaciones de
los par\'ametros de Schur},
Doctoral Dissertation, Universidad de Zaragoza, 1997.

\item{[5]}
T.S. Chihara,
{\sl An Introduction to Orthogonal Polynomials},
Gordon and Breach, New York, 1978.

\item{[6]}
A.J. Dur\'an,
On orthogonal polynomials with respect to a positive definite matrix of measures,
{\sl Can. J. Math.} {\bf 47} (1995) 88--112.

\item{[7]}
A.J. Dur\'an,
Ratio asymptotic for orthogonal matrix polynomials,
{\sl J. Approx. Theory} {\bf 100} (1999) 304--344.

\item{[8]}
A.J. Dur\'an, W. Van Assche,
Orthogonal matrix polynomials and higher-order recurrence relations,
{\sl Linear Alg. Appl.} {\bf 219} (1995) 261--280

\item{[9]}
G. Freud,
{\sl Orthogonal Polynomials},
Akad\'emiai Kiad\'o, Budapest and Pergamon Press, Oxford, 1971, 1985.

\item{[10]}
P.A. Fuhrmann,
Orthogonal matrix polynomials and system theory,
Conference on linear and nonlinear mathematical control theory,
{\sl Rend. Sem. Mat. Univ. Politec. Torino} (1987) 68--124.

\item{[11]}
J. Geronimo,
Scattering theory and matrix orthogonal polynomials on the real line,
{\sl Circuits systems Signal Process} {\bf 1} (1982) 471--494.

\item{[12]}
G.H. Golub, C.F. Van Loan, {\sl Matrix  Computations}, The Johns
Hopkins University Press, Baltimore, second ed., 1989.

\item{[13]}
G.H. Golub, R. Underwood,
The block Lanczos methods for computing eigenvalues,
in: J.R. Rice, Ed., {\sl Mathematical Software III}, Academic Press, New York (1977)
364--377.

\item{[14]}
F.Marcell\'an, G. Sansigre,
On a class of matrix orthogonal polynomials on the real line,
{\sl Linear Alg. Appl.} {\bf 181} (1993) 97--109.

\item{[15]}
F. Marcell\'an, I. Rodr\'{\i}guez,
A class of matrix orthogonal polynomials on the unit circle,
{\sl Linear Alg. Appl.} {\bf 121} (1989) 233--241

\item{[16]}
F. Marcell\'an, A. Ronveaux,
On a class of polynomials orthogonal with respect to a discrete Sobolev inner product,
{\sl Indag. Math. (N.S.)} {\bf 1} (1990) 451--464.

\item{[17]}
F. Marcell\'an, H.O. Yakhlef,
Recent trends on analytic properties of matrix orthonormal polynomials,
{\sl Manuscript}.

\item{[18]}
A. Mat\'e, P. Nevai, V. Totik,
Extensions of Szeg\H o's theory of orthogonal polynomials, II, III,
{\sl Constr. Approx.} {\bf 3} (1987) 51--72; 73--96.

\item{[19]}
E.A. Rakhmanov,
On the asymptotics of the ratio of orthogonal polynomials,
{\sl Math. USSR Sb.} {\bf 32} (1977) 199--213.

\item{[20]}
E.A. Rakhmanov,
On the asymptotics of the ratio of orthogonal polynomials II,
{\sl Math. USSR Sb.} {\bf 46} (1983) 105--117.

\item{[21]}
L. Rodman,
Orthogonal matrix polynomials,
in: P. Nevai, Ed., {\sl Orthogonal polynomials: Theory and Practice}, vol. 294 of
NATO ASI Series C, Kluwer, Dordrecht (1990) 345--362.

\item{[22]}
W. Rudin,
{\sl Real and Complex Analysis},
Mc Graw-Hill, New York, 1966.

\item{[23]}
A. Sinap, W. Van Assche,
Orthogonal matrix polynomials and applications,
Proceedings of the Sixth International Congress on Computational and Applied Mathematics
(Leuven, 1994),
{\sl J. Comput. Appl. Math.} {\bf 66} (1996) 27--52.

\item{[24]}
G. Szeg\H o,
{\sl Orthogonal Polynomials},
4th ed., AMS Colloq. Publ., vol. 23, Amer. Math. Soc., Providence, RI, 1975.

\item{[25]}
H.O. Yakhlef, F. Marcell\'an, M. Pi\~nar,
Relative asymptotics for orthogonal matrix polynomials with convergent recurrence
coefficients,
{\sl J. Approx. Theory}. To appear.

 \bye